\documentclass[11pt,leqno]{amsart}
\usepackage{amsmath,amsfonts,amssymb,amscd,amsthm,amsbsy,upref,}
\numberwithin{equation}{section}
\newtheorem{thm}{Theorem}[section]
\newtheorem{lem}[thm]{Lemma}
\newtheorem{cor}[thm]{Corollary}
\newtheorem{prop}[thm]{Proposition}
\theoremstyle{definition}
\newtheorem{defn}[thm]{Definition}
\newtheorem{prob}[thm]{Problem}             
\newtheorem{remark}[thm]{Remark}
\newtheorem{fact}[thm]{Fact}
\def\nat{{\mathbb N}}
\def\real{{\mathbb R}}
\def\P{{\mathcal P}}
\def\ep{\varepsilon}
\catcode`@=11 \@mparswitchfalse  
\newcounter{mnotecount}[section]

\newcommand{\mnote}[1]{} 

\begin{document}
\title{Lattice Structures and Spreading Models}
\author{S. J. Dilworth, E. Odell \and B. Sar\i}
\address{Department of Mathematics\\ University of South
Carolina, Columbia, SC 29208-0001}
\email{dilworth@math.sc.edu}
\address{Department of Mathematics\\ The University of Texas at Austin\\
1 University Station C1200\\ Austin, TX 78712-0257}
\email{odell@math.utexas.edu}
\address{Department of Mathematics\\ University of North Texas, Denton,
TX 76203-1430} \email{bunyamin@unt.edu}

\thanks{Research of the second
named author was partially supported by the National Science
Foundation. The third named author had a visiting appointment at
the University of South Carolina for the 2004-05 academic year
during part of his research.}

\begin{abstract}
We consider problems concerning the partial order structure of the
set of spreading models of Banach spaces. We construct examples of
spaces showing that the possible structure of these sets include
certain classes of finite semi-lattices and countable lattices, and
all finite lattices.
\end{abstract}
\maketitle

\baselineskip=18pt
\tableofcontents            
\setcounter{section}{-1}
\section{Introduction}

The spreading models of a Banach space $X$ usually have a simpler
and better structure, both individually and collectively, than the
class of subspaces of $X$. Sometimes knowledge of the spreading
models can be used to deduce subspace knowledge about $X$ itself
(e.g., \cite{AOST, OS1}) but the relationship is still not
completely understood. Spaces with no ``nice'' subspaces can have
very nice spreading models (e.g., \cite{AD}).

In this paper we explore further the relationship between a space
and its spreading models. In particular we study the possible
partial order structures of the spreading models of $X$ generated by
normalized weakly null sequences. In \S1 we recall what is known and
unknown and present some new structural observations along with the
relevant background. In \S2 and \S3 we construct spaces $X$ with
certain prescribed spreading model structures. In \S2 we construct
for each $n\in\nat$, a space $X_n$ with $(SP_w(X_n),\le)$ order
isomorphic  to $(\P(n)\setminus\{\emptyset\},\subseteq)$ where
$\P(n)$ is the power set of $\{1,\ldots,n\}$. In \S3 we show that if
$L$ is a countable lattice with a minimum element not containing an
infinite strictly increasing sequence, then there exists a reflexive
space $X_L$ with $(SP_w(X_L),\le)$ order-isomorphic to $L$. The
construction uses some beautiful classical work of Lindenstrauss and
Tzafriri \cite{LT} on Orlicz sequence spaces.

\section{Background, questions and observations}

We use standard Banach space notation and terminology as in \cite{LT}.

Let $X$ be a separable infinite-dimensional Banach space.
A normalized basic sequence $(x_i) \subseteq X$ generates a
{\em spreading model} $(\tilde x_i)$ if for some $\ep_n\downarrow 0$
for all $n\in\nat$ and $(a_i)_1^n \subseteq [-1,1]$,
$$(1+\ep_n)^{-1}\Big\| \sum_{i=1}^n a_i \tilde x_i\Big\|
\le \Big\| \sum_{i=1}^n a_i x_{k_i} \Big\| \le (1+\ep_n) \Big\|
\sum_{i=1}^n a_i \tilde x_i\Big\|$$ for all $n\le k_1<\cdots < k_n$.
$(\tilde x_i)$ is a basic sequence which is 1-spreading and
suppression-$1$ unconditional if $(x_i)$ is weakly null. Every
normalized basic sequence has a subsequence which generates a
spreading model (see \cite{BL} for these and more elementary facts
about spreading models).

We let $[(\tilde x_i)]$ denote the equivalence class of all
spreading models of $X$ which are equivalent (see below) to
$(\tilde x_i)$. $SP_w(X)$ denotes the set of all such $[(\tilde
x_i)]$ where we restrict ourselves only to spreading models
generated by weakly null sequences. If $SP_w(X)=\emptyset$ then
$X$ is  a Schur space, so every normalized spreading model of $X$
is equivalent to the unit vector basis of $\ell_1$ by Rosenthal's
$\ell_1$ theorem \cite{R}.

If $[(\tilde x_i)], [(\tilde y_i)] \in SP_w(X)$ we write $[(\tilde
x_i)] \le [(\tilde y_i)]$ if for some $C<\infty$, $(\tilde y_i)$
\textit{$C$-dominates} $(\tilde x_i)$, i.e., for all $(a_i)\subseteq
\real$ $$\| \sum a_i\tilde x_i\| \le C\| \sum a_i \tilde y_i\|\ .$$
$(\tilde x_i)$ and $(\tilde y_i)$ are \textit{equivalent} if each
dominates the other. $(SP_w(X),\le)$ is a partially ordered set.

We sometimes have occasion to consider a specific $(\tilde x_i)$ and
shall abuse notation by writing ``let $(\tilde x_i)\in SP_w(X)$.''
$c_{00}$ denotes the linear space of finitely supported real
sequences.

\begin{fact} \label{fact1.1}    \cite{AOST}
$(SP_w(X),\le)$ is a \textit{semi-lattice}, i.e., each two elements of
$SP_w(X)$ admit a least upper bound. Moreover if $(\tilde
x_i),(\tilde y_i) \in SP_w(X)$ there exists $(\tilde z_i)\in
SP_w(X)$ which is 2-equivalent to the subsymmetric norm on $c_{00}$
given by
$$\|(a_i)\| = \|\sum a_i \tilde x_i\| \vee \|\sum a_i\tilde y_i\|.$$
\end{fact}
\begin{fact} \label{fact1.2}\cite{AOST}        
Every countable subset of $(SP_w(X),\le)$ admits an upper bound.
Moreover if $(\tilde x_i^n)_{i=1}^\infty \in SP_w(X)$ for $n\in\nat$
and $(C_n)_{n=1}^\infty \subseteq (0,1)$ with $\sum_{n=1}^\infty
C_n^{-1}\le1$ then there exists $(\tilde z_i) \in SP_w(X)$ which
$C_n$-dominates $(\tilde x_i^n)$ for each $n\in\nat$. In addition
for $(a_i)\in c_{00}$
$$\|\sum a_i \tilde z_i\| \le
C_1^{-1} \bigg(\sum_{n=1}^\infty C_n^{-1} \|\sum a_i \tilde x_i^n\|
\bigg)\ .$$ \ We shall designate this $(\tilde z_i)$ by the notation
$(\tilde z_i) = (\sum C_n^{-1} \tilde x_i^n)$, which in fact is
motivated by the proof in \cite{AOST} (the precise quantification as
given above is noted in \cite{S2}). \end{fact}
\begin{fact} \label{fact1.3} \cite{S2}          
If $SP_w(X)$  admits an infinite strictly increasing sequence then
$SP_w(X)$ is uncountable. In fact there exist $[(\tilde y_i^\alpha)]
\in SP_w(X)$ for $\alpha<\omega_1$ so that  $[(\tilde y_i^\alpha)] <
[(\tilde y_i^\beta)]$ if $\alpha < \beta <\omega_1$. \end{fact}

Our next result is motivated by the proof of Fact~\ref{fact1.3}.

\begin{thm}    \label{thm1.4}
Let $I$ be an infinite set and
let $(\tilde x_i^\alpha)_{i=1}^\infty \in SP_w(X)$ for $\alpha\in
I$. For $A\subseteq I$ define a subsymmetric norm on $c_{00}$ by
$R_A (a_i) = \sup_{\alpha \in A} \|\sum a_i \tilde x_i^\alpha\|$. If
for every non-empty  finite  $F\subseteq I$,
$R_I$ is not equivalent to $R_F$, then
$SP_w(X)$ admits an infinite strictly increasing sequence.
\end{thm}
\begin{proof}
We may assume $I=\nat$.
We shall construct a strictly increasing sequence
$(\tilde y_i^n)_{i=1}^\infty$ for $n\in\nat$.
We shall let
$$(\tilde y_i^1) = \bigg( \sum_{n=1}^\infty  2^{-n} \tilde
z_i^n\bigg) \qquad \text{(see Fact~\ref{fact1.2})}$$ where $(\tilde
z^n)_{n=1}^\infty$ is a reordering of $(\tilde x^n)_{n=1}^\infty$
selected as follows. Let $\ep_n\downarrow 0$ and for each
$n\in\nat$, $\tilde z^{2n}$ is chosen so that for some
$(a_\ell^n)_{\ell=1}^\infty \in c_{00}$,
$$R_I (a_\ell^n)_{\ell=1}^\infty =1\ ,\quad
\Big\| \sum_\ell a_\ell^n \tilde z_\ell^{2n}\Big\| > \frac12$$
and $R_{I_n} (a_\ell^n)_{\ell=1}^\infty < \ep_n 2^{-2n}$ where
$$I_n = \{ m\in\nat : \tilde x^m = \tilde z^j\ \text{ for some }\
j\le 2n-1\} .$$ $\tilde z^n$ for $n$ odd is selected arbitrarily so
as to exhaust the collection $(\tilde x^s)_{s\in\nat}$.

For $n\in\nat$ we have  (see Fact~\ref{fact1.2})
\begin{equation*}
\begin{split}
\Big\| \sum_\ell a_\ell^n \tilde y_\ell^1\Big\|
& \le 2R_{I_n} (a_\ell^n)_{\ell=1}^\infty
+ 2\cdot 2^{-2n} \Big\|\sum_\ell a_\ell^n \tilde z_\ell^{2n}\Big\|
+ 2\sum_{m>2n} 2^{-m} \Big\| \sum_\ell a_\ell^n \tilde
z_\ell^m\Big\|\\
&< 2\ep_n 2^{-2n} + 2\cdot 2^{-2n} + 2\sum_{m>2n} 2^{-m}\\
& = (2\ep_n +4 ) 2^{-2n}\ .
\end{split}
\end{equation*}
Furthermore
$$\| \sum a_\ell^n \tilde y_\ell^1\|
\ge 2^{-2n} \Big\| \sum_\ell a_\ell^n \tilde z_\ell^{2n}\Big\|
> \frac12 2^{-2n}\ .$$
We thus obtain that $(\tilde x_i^n) < (\tilde y_i^1)$ for all
$n\in\nat$ and $(\tilde y_i^1) < R_I$. Moreover
we can
iterate the argument beginning anew with the collection $\{ (\tilde
y_i^1)\} \cup \{ (\tilde x_i^n)\}_{n=1}^\infty$,
which satisfies the same hypothesis as $(\tilde x_i^n)\}_{n=1}^\infty$,
to obtain $\tilde
y^2$, and so on.
\end{proof}

\begin{fact}    \label{fact1.5}       
$SP_w(X)$ can be hereditarily uncountable \cite{AOST},  i.e.,
$SP_w(Y)$ is uncountable for all infinite-dimensional subspaces $Y$
of $X$. If $SP_w(X)$ is countable, then by a diagonal argument one
can find $X_0 \subseteq X$ with  $SP_w(X_0) = SP_w(Y)$ for all
$Y\subseteq X_0$. It may be that then $|SP_w(X_0)| =1$ but this
remains open.  \end{fact}

We also have the

\begin{prob} \label{prob1.6}
If $X$ is reflexive and $SP_w(X)$ is countable must some $(\tilde
x_i)  \in SP_w(X)$ be equivalent to the unit vector basis of $c_0$
or $\ell_p$ for some $1\le p<\infty$?
\end{prob}

If so this would be a case where one would have a stronger theorem
than Krivine's \cite{K}. Not every reflexive space has a spreading
model isomorphic to $c_0$ or some $\ell_p$ (\cite{OS2},
\cite{AOST}). In the nonreflexive case it is possible to have
$|SP_w(X)| =1$ yet the unique spreading model is not $c_0$ or any
$\ell_p$. This is the case for certain Lorentz sequence spaces
$d_{w,1}$ (see \S2).

Problem~\ref{prob1.6} was raised and partially solved in the case
$|SP_w(X)|=1$ in \cite{AOST}. We give some further partial results
below.
\begin{remark}    \label{remark1.7}     
Assume that $SP_w(X)$ is countable or more generally does not admit
an infinite strictly increasing sequence. For $[(\tilde x_i)]\in
SP_w(X)$ and $(a_i) \in c_{00}$, define
$$R(a_i) \equiv R_{[(\tilde x_i)]} (a_i)
= \sup \left\{ \Big\| \sum_1^n a_i\tilde y_i\Big\| : (\tilde y_i)
\in [(\tilde x_i)]\right\}\ .$$ By 1.4  $R$ is equivalent to
$(\tilde x_i)$. Thus for each $(\tilde y_i) \in [(\tilde x_i)]$
there exists $C<\infty$ so that $(\tilde y_i)$ $C$-dominates every
$(\tilde z_i)\in [(\tilde x_i)]$. Also \cite{S2} there exists $p= p
(\tilde x_i) \in [1,\infty]$ so that for all $1\le q < p$ there
exists $C_q<\infty$ so that for all $(a_i)\subseteq \real$,
$$\bigg(\sum |a_i|^p\bigg)^{1/p}
\le R(a_i) \le C_q\bigg(\sum |a_i|^q\bigg)^{1/q}\ .$$ $p(\tilde
x_i)$ is the infimum of the ``Krivine $p$'s'' for $(\tilde x_i)$
(see \cite{S2}). It is mistakenly stated in \cite{S2} that, in this
case, $p(\tilde x_i)$ is the only Krivine $p$. However, this is not
yet clear. \end{remark}

\begin{remark}   \label{remark1.8}     
Let $SP_w(X)$ be stabilized hereditarily for $X$. Then
for all $(\tilde x_i)\in SP_w(X)$ there exist
$X_0 \subseteq X$ and
$C<\infty$ such that: for all $Y\subseteq X_0$ there exists
$(\tilde y_i) \in SP_w(Y)$ which is $C$-equivalent to
$(\tilde x_i)$.

The proof is elementary.
Assume not and use a diagonal argument to get a contradiction.
\end{remark}

\begin{thm} \label{thm1.9}
Suppose that $SP_w(X)$ is  countable and that $SP_w(Y^*)$
is countable for all infinite-dimensional subspaces $Y$ of $X$. Then every
$(\tilde e_i) \in SP_w(X)$
is  equivalent to the unit vector basis of $c_0$ or $\ell_p$
for some $1 \le p < \infty$.\end{thm}
\begin{proof} Let $(e_i)$ be a normalized weakly null sequence in $X$
generating the spreading model $(\tilde e_i)$. By passing to a
subsequence and renorming we may assume that $(e_i)$ is bimonotone
basic and Schreier-unconditional, i.e.\ for some $\ep_n\downarrow 0$
and all $F\in S_1$ and $(a_i)\in c_{00}$, \begin{equation}
\label{eq: schreieruncond} \Big\| \sum_{i\in F} a_i e_i\Big\| \le
(1+\ep_{\min F}) \| \sum a_i e_i\| \end{equation} (see \cite{O}).
Here $F\in S_1$ (first Schreier class) if $|F|\le \min F$.

We may assume that no subsequence of $(e_i)$ is equivalent to the
unit vector basis of $c_0$. Thus by passing to a further subsequence
we may assume that  $(f_i)$, the sequence of biorthogonal functions
to $(e_i)$, is weakly null in $[(e_i)]^*$ \cite[Cor. 4.4]{O2}. From
\eqref{eq:  schreieruncond} it is easy see that $(f_i)$ is
normalized and has spreading model $(\tilde f_i)$ which is
$1$-equivalent to $(\tilde e_i^*)$, the biorthogonal functionals to
$(\tilde e_i)$ in $[(\tilde e_i^*)]$. By Krivine's theorem \cite{K}
it suffices to prove that, for some $D<\infty$, every spreading
model $(\tilde x_i)$ of an identically distributed block basis
$(x_i)$ of $(e_i)$ with $\|x_i\| \rightarrow 1$ is $D$-equivalent to
$(\tilde e_i)$. Note that $(x_i)$ is weakly null and $(\tilde x_i)$
is equivalent to an identically distributed normalized block basis
of $(\tilde e_i)$ and hence to $(\tilde e_i)$. Since $SP_w(X)$ is
countable, by Theorem~\ref{thm1.4} there exists $C_1<\infty$ (which
depends only on $\tilde e_i)$) such that
\begin{equation} \label{eq: upest}
\|\sum a_i \tilde x_i\| \le C_1\|\sum a_i \tilde e_i\| \qquad ((a_i)
\in c_{00}). \end{equation} We may choose an identically distributed
block basis $(g_i)$ of $(f_i)$ with $\operatorname{supp} (g_i)
\subseteq \operatorname{supp} (x_i)$, $\|g_i\| \rightarrow 1$, and
$g_i(x_i) \rightarrow 1$. Note that $(g_i)$ has spreading model
$(\tilde g_i)$ which is $1$-equivalent to an identically distributed
block basis of $(\tilde f_i)$. Also $(g_i)$ is weakly null and since
$SP_w(X^*)$ is countable we have, again by  Theorem~\ref{thm1.4},
that there exists $C_2<\infty$ (which depends only on $(\tilde
f_i)$) such that
$$\|\sum a_i \tilde g_i\| \le C_2\|\sum a_i \tilde f_i\| \qquad
((a_i) \in c_{00}).$$ Let $h_i$ be the restriction of $g_i$ to
$[(x_i)]$. Since $(x_i)$ is bimonotone and Schreier unconditional,
we have as above that $(h_i)$ has spreading model $(\tilde h_i)$ in
$[(x_i)]^*$ which is $1$-equivalent to $(\tilde x_i^*)$, the
biorthogonal functionals to $(\tilde x_i)$. Thus for $(a_i) \in
c_{00}$,
$$ \|\sum a_i \tilde x_i^*\|=\|\sum a_i \tilde h_i\|
\le \|\sum a_i \tilde g_i\| \le C_2 \|\sum a_i \tilde f_i\|.$$ By
duality, \begin{equation} \label{eq: lowest} \|\sum a_i \tilde x_i\|
\ge  \frac{1}{C_2}\|\sum a_i \tilde e_i\| \qquad ((a_i) \in c_{00}).
\end{equation}
Thus by \eqref{eq: upest} and \eqref{eq: lowest},  $(\tilde x_i)$ is
$D \equiv C_1C_2$-equivalent to $(\tilde e_i)$.
\end{proof}

\begin{thm}   \label{thm1.10}      
Let $X$ be reflexive with  $|SP_w(X)| = |SP_w(X^*)| =1$. Assume also
that the element of $SP_w(X^*)$ is equivalent to the biorthogonal
functionals of the element $(\tilde x_i)$ in $SP_w(X)$. Then
$(\tilde x_i)$ is equivalent to the unit vector basis of $c_0$ or
$\ell_p$ for some $1\le p<\infty$.
\end{thm}
\begin{proof}
We first note that if $X_0$ is any infinite-dimensional subspace of
$X$ then $X_0$ satisfies the same hypothesis as $X$. Indeed the only
question here is the uniqueness of the spreading  models in $X_0^*$.
Let $(\tilde f_i)$ be a normalized spreading model for $X_0^*$
generated by $(f_i)$. Then $(f_i)$ is the image under the quotient
map of a seminormalized weakly null sequence in $X^*$ and this
yields that $(\tilde x_i^*)$ dominates $(\tilde f_i)$. A similar
argument applied to  the sequence biorthogonal to
$(f_i)$ shows that
$(\tilde f_i)$  dominates $(\tilde
x_i^*)$. The result now follows from Theorem~\ref{thm1.9}.
\end{proof}

The proof of Theorem~\ref{thm1.9} contains the following result.

\begin{thm}  \label{thm1.11}       
Let $(e_i)$ be a normalized basis for a reflexive space $X$ which is
$C$-Schreier unconditional for some $C<\infty$, i.e.,
$$\Big\| \sum_F a_i e_i\Big\|
\le C\|\sum a_i e_i\|\ \text{ for all }\ F\in S_1$$ and
$(a_i)\subseteq \real$. If $|SP_w(X)| = |SP_w(X^*)| =1$ then the
unique spreading model of $X$ is equivalent to the unit vector basis
of $c_0$ or $\ell_p$ for some $1\le p<\infty$.
\end{thm}
\begin{remark}  \label{remark1.12}         
If $SP_w(X)$ is countably infinite then $SP_w(X)$ contains
$\{(\tilde x_i^n)_{i=1}^\infty :n\in\nat\}$ with either $(\tilde
x_i^n) > (\tilde x_i^m)$ for all $n<m$ or $(\tilde x_i^n)$ and
$(\tilde x_i^m)$  mutually incomparable for all $n\ne m$. Indeed
Ramsey's theorem yields a subsequence of any sequence of spreading
models satisfying either one of the two possibilities above or a
sequence that is strictly increasing. The latter is ruled out by
Fact~\ref{fact1.3}. Both possibilities can occur for reflexive
spaces. As noted elsewhere \cite{AOST} (see also Theorem~\ref{prop:
ordinal} below) it is easy to check that every spreading model of
$(\sum \oplus \ell_{p_n})_{p_1}$ is equivalent to some $\ell_{p_n}$
if $p_1 < p_2 < \cdots$. In \S3 we shall  show the second (mutually
incomparable) possibility.

The uncountable case is less clear.
\end{remark}
\begin{prob}   \label{prob1.13}         
If $SP_w(X)$ is uncountable must there exist $\{(\tilde
x_i^\alpha)_{i=1}^\infty : \alpha <\omega_1\} \subseteq SP_w(X)$
which is either strictly increasing w.r.t.\ $\alpha$, strictly
decreasing or consists of mutually incomparable elements.
\end{prob}

If there is a counterexample, $X$, say, to this question, then by
Fact~\ref{fact1.3} $SP_w(X)$ cannot contain an infinite increasing
sequence. We do not know however the answer to this generalized
version of Problem~\ref{prob1.13}.

\begin{prob} \label{newprob1.14}
Let $L$ be an uncountable semi-lattice which admits no infinite
strictly increasing sequence.
Must $L$ admit a family $(x_\alpha)_{\alpha<\omega_1}$ with either
\begin{itemize}
\item[(i)] $\forall\ \alpha\ne\beta$, $x_\alpha$ and $x_\beta$ are
incomparable or
\item[(ii)] $\forall\ \alpha <\beta<\omega_1$, $x_\alpha >x_\beta$?
\end{itemize}
\end{prob}

The following example due to Sierpinski (see \cite{ER}) provides a
counterexample to the corresponding question for posets. Let $L =
(\{x_\alpha\}_{\alpha<\omega_1},\preceq)$, where
$\{x_\alpha\}_{\alpha<\omega_1}$ are distinct points in $(0,1)$ with
$x_\alpha \prec x_\beta$ iff $\alpha >\beta$ and $x_\alpha
> x_\beta$ (in $\real$). Then $L$ is a poset without any infinite
increasing sequences and without any uncountable chains or
antichains.

If an $\omega_1$-Suslin tree exists then Problem~\ref{newprob1.14}
easily has a negative answer.
 In a related result Shelah
\cite{Sh}  has shown that under (CH) there exists an uncountable
Boolean algebra without uncountable chains or antichains and
moreover, (CH) + no $\omega_1$-Suslin tree is consistent with ZFC.
In particular, under (CH) there is a counterexample to
Problem~\ref{prob1.13}  if $SP_w(X)$ is replaced by a general
semi-lattice.

 Our work in the next two sections suggests
the following.

\begin{prob} \label{prob1.14}           
Let $L$ be a countable semi-lattice not admitting an infinite
strictly increasing sequence. Does there exist $X$ (possibly even
reflexive) with $(SP_w(X),\le)$ order-isomorphic to $L$? \end{prob}

\section{Spreading model sets without a minimum element}

In this section we shall construct some families of Banach spaces
whose spreading model sets  do not have a minimum element in the
domination ordering. The Banach spaces in question are finite direct
sums of certain Lorentz sequence spaces $d(w,p)$.

The construction depends on the existence of an arbitrary number
of incomparable submultiplicative functions. We begin with a
technical  definition to facilitate the discussion.
\begin{defn}\label{def: submultiplicative}
Let $2 \le n_0 \le \infty$ and let $S$ be a real-valued function
defined on $[1,n_0]$. We shall say that $S$ is
\textit{submultiplicative on $[1,n_0]$} (or on $[1,\infty)$ if
$n_0=\infty$) if $S$ satisfies the following conditions:
\begin{itemize}
\item[(a)] $S$ is piecewise-linear, continuous, strictly
increasing, and concave. \item[(b)] $S(x)=x$ for $1 \le x \le 2$.
\item[(c)]$S(xy) \le S(x)S(y)$ for all $x,y$ such that $1 \le x,y,
xy \le n_0$.
\end{itemize} \end{defn}
\begin{lem} \label{lem: slowdown}  Suppose that $2\le n_0  < \infty$
and that $S$ is submultiplicative on $[1,n_0]$. Then there exists
$\varepsilon_0 >0$ such that for all $0 < \varepsilon <
\varepsilon_0$ the extension $S_\varepsilon$ of $S$ to the interval
$[1, n_0^2]$ defined by \begin{equation*}
S_\varepsilon(x) = \begin{cases} S(x) &\text{for $1\le x \le n_0$}\\
S(n_0) + \varepsilon(x-n_0) &\text{for $n_0<x\le n_0^2$} \end{cases}
\end{equation*}
is submultiplicative on $[1,n_0^2]$.
\end{lem}
\begin{proof} Since $S$ is continuous, piecewise-linear, and
strictly increasing on $[1,n_0]$, there exists $c>0$ such that
\begin{equation} \label{eq: sublinear} S(x) \ge S(x-h) + ch \qquad
(1 \le x-h \le x  \le n_0).
\end{equation} Define $\tilde S$ on $[1,n_0^2]$ as follows:
\begin{equation*}
\tilde S(x) := \inf \{S(a)S(b) \colon x=ab, 1 \le a,b\le n_0\}.
\end{equation*} Since $S$ is continuous and strictly increasing
on $[1,n_0]$ it follows that $\tilde S$ is continuous and strictly
increasing on $[1,n_0^2]$. Moreover, conditions (b) and (c) of
Definition~\ref{def: submultiplicative} imply that $\tilde S(x) =
S(x)$ for all $x \in [1,n_0]$. Suppose that $x := n_0+h$ satisfies
$n_0 \le x \le n_0^2$. By compactness there exist $a_x,b_x$ such
that $\tilde S(x) = S(a_x)S(b_x)$, $x=a_xb_x$, and $1 \le a_x \le
b_x \le n_0$. Then
\begin{align*}
S(n_0) &= \tilde S(n_0)
\le S(a_x) S(b_x - \frac{h}{a_x})\\ \intertext{(since
$n_0=a_x(b_x-h/a_x)$ and $S$ is submultiplicative on $[1,n_0]$)}
&\le S(a_x)(S(b_x)- \frac{ch}{a_x})\\
\intertext{(by \eqref{eq: sublinear})}
&=\tilde S(n_0+h)- \frac{cS(a_x)}{a_x} h\\
&\le \tilde S(n_0+h)- \frac{c}{n_0} h,
\end{align*} (since $S(a_x)\ge1$ and $a_x \le n_0$).
So $$S(n_0) + \frac{c}{n_0}h \le \tilde S(n_0+h).$$ Hence,
provided $\varepsilon < c/n_0$, we have $S_\varepsilon(x) \le
\tilde S(x)$ for $1 \le x \le n_0^2$. To verify
submultiplicativity of $S_\varepsilon$ on the interval $[1,n_0^2]$
it remains to check that $$S_\varepsilon(xy) \le
S_\varepsilon(x)S_\varepsilon(y)$$ for all $1 \le x \le n_o$ and
$n_0 \le y \le n_0^2$ such that $xy \le n_0^2$. Since
$S_\varepsilon(xy) = S_\varepsilon(y) + \varepsilon(xy-y)$, we
require
\begin{equation*}
\frac{S_\varepsilon(y)+\varepsilon(x-1)y}{S_\varepsilon(y)} \le
S_\varepsilon(x),
\end{equation*} i.e. \begin{equation} \label{eq: epsiloncondition}
\varepsilon(x-1)y \le (S_\varepsilon(x)-1)S_\varepsilon(y).
\end{equation} First consider the case $x \ge 2$. Then
$S_\varepsilon(x)-1\ge1$ and since $(x-1)y\le xy \le n_0^2$ it
follows that  \eqref{eq: epsiloncondition} will be satisfied
provided $n_0^2\varepsilon \le S(n_0)$.  On the other hand, if $1
\le x \le 2$, then by condition (b) of Definition \ref{def:
submultiplicative}, \eqref{eq: epsiloncondition} reduces to
$\varepsilon y \le S_\varepsilon(y)$, which will again be satisfied
provided $n_0^2\varepsilon \le S(n_0)$. This proves the lemma for
$$ \varepsilon_0 = \min(\frac{c}{n_0},\frac{S(n_0)}{n_0^2}).$$
\end{proof}

By an obvious repeated application of Lemma~\ref{lem: slowdown} one
obtains the following result.
\begin{lem} \label{lem: slowdown2} Suppose that $n_0\ge2$ and that
  $S$ is submultiplicative on $[1,n_0]$.
Then, given $\varepsilon>0$ and $N_0>n_0$, there exists a
submultiplicative extension of
$S$ to $[1,N_0]$ such that $S(N_0) < S(n_0)+\varepsilon$. \end{lem}
\begin{lem} \label{lem: speedup}
Suppose that $S$ is submultiplicative on $[1,n_0]$, where $n_0\ge2$,
and that $S(n_0) = K \ge 2$. Then there exist $N_0 > n_0$ and a
submultiplicative extension of $S$ to $[1,N_0]$ such that $S(N_0)
\ge 3K/2$. \end{lem}
\begin{proof} Let $n_1=n_0^2$. By Lemma~\ref{lem: slowdown} we may and
shall assume that $S$ has been extended to be submultiplicative on
$[1,n_1]$. By a second application of Lemma~\ref{lem: slowdown}
there exists $\varepsilon>0$ such that
\begin{equation*} S_\varepsilon(x) = \begin{cases} S(x) &\text{for $1\le
      x \le n_1$}\\
S(n_1) + \varepsilon(x-n_1) &\text{for $n_1<x\le n_1^2$} \end{cases}
\end{equation*}
is submultiplicative on $[1,2n_1]$ (or even $[1,n_1^2]$ although we
will only use submultiplicativity on $[1,2n_1]$). If
$S_\varepsilon(2n_1) \ge 3K/2$ then we are done.
So we may assume that $S(2n_1) < 3K/2$, which implies (since $S(n_1)
\ge K$) that
\begin{equation} \label{eq: epsiloncondition2}
n_1 \varepsilon< \frac{K}{2}. \end{equation} Choose $N_0>2n_1$ such
that $S_\varepsilon(N_0) = 3K/2$. We shall show that $S_\varepsilon$
is submultiplicative on $[1,N_0]$. So suppose that $1 \le x \le y
\le xy \le N_0$. Since $S_\varepsilon(x)$ is submultiplicative on
$[1,2n_1]$ we may assume that $xy \ge 2n_1$. Since $n_1 = n_0^2$, it
follows that $y \ge n_0$, so $S_\varepsilon(y) \ge K$. First
consider the case $x>2$. Then $S_\varepsilon(x)\ge2$, so
$$S_\varepsilon(xy)\le S_\varepsilon(N_0) < 2K \le
S_\varepsilon(x)S_\varepsilon(y).$$ On the other hand, if $1 \le x
\le 2$, then by condition (b) of Definition~\ref{def:
submultiplicative} $S_\varepsilon(x)=x$,  so the submultiplicativity
condition becomes $$S_\varepsilon(y)+\varepsilon(x-1)y =
S_\varepsilon(xy) \le S_\varepsilon(x)S_\varepsilon(y)
=xS_\varepsilon(y),$$ i.e. $\varepsilon y \le S_\varepsilon(y)$.
This is clearly satisfied if $n_0 \le y \le K/\varepsilon$ since
$S_\varepsilon(y) \ge S_\varepsilon(n_0)=K.$ But \begin{align*}
S_\varepsilon(\frac{K}{\varepsilon}) &= S_\varepsilon(n_1) +
\varepsilon(\frac{K}{\varepsilon}-n_1)\\
&\ge K + K - n_1\varepsilon
\ge K + K - \frac{K}{2}\\ &= \frac{3K}{2}= S_\varepsilon(N_0),
\end{align*}
where the last inequality follows from \eqref{eq: epsiloncondition2}.
Thus, $K/\varepsilon \ge N_0$, which proves that $S_\varepsilon$ is
submultiplicative on $[1, N_0]$ as desired.
\end{proof}
By an obvious repeated application of Lemma~\ref{lem: speedup} one
obtains the following result.
\begin{lem} \label{lem: speedup2} Suppose that $n_0 \ge 2$ and that
$S$ is submultiplicative on $[1,n_0]$. Then, given $M>0$, there
exist $N_0>n_0$ and a submultiplicative extension of $S$ to
$[1,N_0]$ such that $S(N_0)>M$. \end{lem} Next we construct an
infinite collection of mutually incomparable submultiplicative
functions. This will be used to construct spreading model diagrams
in Theorems~\ref{thm: nonreflexivecase} and \ref{thm: reflexivecase}
below. (In fact, the existence of arbitrarily large finite sets of
incomparable submultiplicative functions would suffice for the
applications.)
\begin{prop} \label{prop: constructsubmult} There exists a sequence
  $(S_i)_{i=1}^\infty$ of  submultiplicative functions  on
  $[1,\infty)$ such that for \textit{every}  nonempty finite set
$A\subset\mathbb{N}$ and for every $j \in \mathbb{N}\setminus A$, we
have
\begin{equation} \label{eq: incomparable} \sup_{n \ge 1}
  \frac{S_j(n)}{\max_{i \in A} S_i(n)} = \infty. \end{equation}
\end{prop}
\begin{proof} We shall define $(S_i)_{i=1}^\infty$ on $[1,\infty)$
by defining their values inductively on an increasing sequence of
initial segments $[1,n_0]$. Let us describe the inductive step.
Suppose that $(S_i)$ have  been defined to be submultiplicative on
some initial segment $[1,n_0]$ in such a way that the collection of
restrictions of $(S_i)$ to $[1,n_0]$ is  a \textit{finite} collection
of functions on $[1,n_0]$. Now fix a finite set $A \subseteq
\mathbb{N}$ and  a positive integer $N$. By applying Lemma~\ref{lem:
slowdown2} to $S_i$ ($i\in A$) on $[1,n_0]$, and applying
Lemma~\ref{lem: speedup2} to $S_i$ ($i \in \mathbb{N}\setminus A$)
on $[1,n_0]$ (which is a finite collection by assumption), there
exist $N_0>n_0$ and submultiplicative extensions of $S_i$ to
$[1,N_0]$ such that
$$ \max_{n \in [n_0,N_0]} \frac{S_j(n)}{\max_{i \in A}S_i(n)}>N$$
for all  $j \in \mathbb{N} \setminus A$. At the end of the inductive
step we have defined $(S_i)$ to be submultiplicative on $[1,N_0]$.
Moreover, the new collection of initial segments of
$(S_i)_{i=1}^\infty$ on $[1,N_0]$ thus obtained will be finite. Now
one simply enumerates (in any manner) the countable collection of
possible choices for $A$ and $N$ to carry out the inductive
definition.
\end{proof}

Let $1 \le p < \infty$ and let $w = (w(n))_{n=1}^\infty$ be a
non-increasing sequence of positive \textit{weights} such that
$w(1)=1$, $w(n) \rightarrow 0$ as $n \rightarrow \infty$, and
$\sum_{n=1}^\infty w(n) = \infty$. Recall that the Lorentz sequence
space $d(w,p)$ is the Banach space with Schauder basis $(e_n)$ whose
norm is defined by
\begin{equation*} \|\sum_{n=1}^\infty a_n e_n\|_{w,p} :=
(\sum_{n=1}^\infty a^{*p}_n w(n))^{1/p},
\end{equation*}
where $(a^*_n)_{n=1}^\infty$ is the nonincreasing rearrangement of
any  scalar sequence $(|a_n|)_{n=1}^\infty$ which converges to zero.
Note that \begin{equation*} \|\sum_{n=1}^\infty a_n e_n\|_{w,p} \le
\|\sum_{n=1}^\infty a_n e_n\|_p := (\sum_{n=1} ^\infty
|a_n|^p)^{1/p}.\end{equation*} The corresponding \textit{fundamental
function} $(S(n))_{n=1}^\infty$ is defined by
$$S(n)  =\| \sum_{i=1}^n e_i\|_{w,p}^p = \sum_{i=1}^n w(i).$$  It is
known that $d(w,p)$ contains  subspaces that are almost isometric
to $\ell_p$ and is reflexive if and only if $1<p<\infty$.

The weight $w$ is said to be \textit{submultiplicative} if there
exists a constant $C$ such that $S(mn) \le C S(m)S(n)$
for all $m,n \in \mathbb{N}$. We require the following theorem due  to
Altshuler,
Casazza and Lin.
\begin{fact} \cite{ACL}
Suppose that $w$ is submultiplicative. Then every
normalized block basis  in $d(w,p)$ has a subsequence which is
equivalent to the unit vector basis of $\ell_p$ or to the unit
vector basis of $d(w,p)$.
\end{fact}
The above theorem has the following immediate corollary.
\begin{cor} \label{cor: spreadsubmult}
Suppose that $w$ is submultiplicative. Then every spreading model
of $d(w,1)$ generated by a weakly null sequence is equivalent to
the unit vector basis of $d(w,1)$. For $1<p<\infty$, every
spreading model of $d(w,p)$ generated by a weakly null sequence is
equivalent to the unit vector basis of $\ell_p$ or to the unit
vector basis of $d(w,p)$.
\end{cor}

Note that to each submultiplicative function $S$ defined on
$[1,\infty)$ there corresponds a submultiplicative weight sequence
$w(n) := S(n)-S(n-1)$ (with constant $C=1$) whose fundamental
function is $(S(n))_{n=1}^\infty$. Let $w_i$ ($1 \le i < \infty$) be
the weight sequences corresponding to the submultiplicative
functions $S_i$ constructed in Proposition~\ref{prop:
constructsubmult}. Note that $\lim_{n\rightarrow \infty} w_i(n)=0$
for each $i$.

Now we come to the main results of this section. For $n\in
\mathbb{N}$, let $P(n)$ denote the
power set of $\{1,\dots,n\}$ partially ordered by inclusion.

\begin{thm} \label{thm: nonreflexivecase} For each $n \in
\mathbb{N}$, let $X_n(1):=(\sum_{i=1}^n \oplus d(w_i,1))_\infty$.
Then $SP_w(X_n(1))$ is order-isomorphic to $P(n) \setminus
\{\emptyset\}$. \end{thm}
\begin{proof} Let $(f_j)_{j=1}^\infty$ be a normalized  spreading
model for $X_n(1)$ generated by a weakly null sequence. Then there
exist a nonempty $A \subseteq \{1,\dots,n\}$ and normalized
spreading models $(f^i_j)_{j=1}^\infty$ of $d(w_i,1)$ ($i \in A$),
generated by weakly null sequences in $d(w_i,1)$, such that
\begin{equation*} \|\sum_{j=1}^\infty a_j f_j\| \approx  \max_{i \in
A} \|\sum_{j=1}^\infty a_j f^i_j\|.
\end{equation*}
Thus, by the first part of  Corollary~\ref{cor: spreadsubmult},
\begin{equation} \label{eq: spreadtype}
\|\sum_{j=1}^\infty a_j f_j\| \approx  \max_{i \in A}
\|\sum_{j=1}^\infty a_j e_j\|_{w_i,1}.
\end{equation}
Conversely, the right-hand side of \eqref{eq: spreadtype} defines a
normalized spreading model $SP(A)$
for every nonempty $A \subseteq \{1,\dots,n\}$. Note that
\begin{equation*}
\|\sum_{j=1}^m  f_j\| \approx  \max_{i \in A} S_i(m) \qquad(m \in
\mathbb{N}).
\end{equation*}
Thus, by \eqref{eq: incomparable} of Proposition~\ref{prop:
  constructsubmult}, we have
$$ A \subset B \Leftrightarrow SP(A) < SP(B)$$
for all nonempty $A,B \subseteq \{1,\dots,n\}$.
\end{proof}

In the reflexive case ($1<p<\infty$) we have to add an extra node on
the top.

\begin{thm} \label{thm: reflexivecase} Let $1<p<\infty$ and, for
each $n \in \mathbb{N}$, let $X_n(p):= (\sum_{i=1}^n \oplus
d(w_i,p))_\infty$. Then $SP_w(X_n(p))$ is order-isomorphic to $(P(n)
\cup \{\{1,\dots,n+1\}\}) \setminus \{\emptyset\}$.
\end{thm}
\begin{proof} The proof is essentially the same as before. However,
from the second part of Corollary~\ref{cor: spreadsubmult}, we
obtain an extra spreading model equivalent to the unit vector basis
of $\ell_p$ which dominates every other spreading model. This
spreading model corresponds to $\{1,\dots,n+1\}$ under the
order-isomorphism.
\end{proof}

\section{Countable lattices with a minimum element}
 Recall that a
\textit{lattice} is a partially ordered set in which any two
elements have both a least upper bound and a greatest lower bound.
The following theorem is the main result of this section.
\begin{thm} \label{thm: lattice}
Let $L$ be a countable lattice with a minimum element not containing
an infinite increasing sequence. Then there exists a reflexive space
$X_L$ such that $SP_w(X_L)$ is order-isomorphic to $L$.
\end{thm}

\begin{remark} Recall that  $(SP_w(X),\le)$ is always a
semi-lattice, i.e. every two elements have a least upper bound
(Fact~\ref{fact1.1}), and that when countable it does not contain
any infinite increasing sequences (Fact~\ref{fact1.3}). It is easy
to see that  such a semi-lattice with a minimum element is
automatically a lattice. Thus, Theorem~\ref{thm: lattice}
characterizes the possible  poset structure of $(SP_w(X), \le)$
when $SP_w(X)$  is countable and has a minimum element.
\end{remark}
The space $X_L$ will be an $\ell_p$ direct-sum of suitably
constructed Orlicz sequence spaces. The proof of the theorem will be
given at the end of the section. First we recall some preliminary
facts about Orlicz spaces. All the unexplained terms and facts can
be found in Chapter 4 of \cite{LT}, with which our notation is
consistent.

An \textit{Orlicz function} $M$ is a real-valued continuous
non-decreasing and convex function defined on $[0,1]$ such that
$M(0)=0$ and $M(1)=1$. For a given $M$, the \textit{Orlicz sequence
space} $\ell_M$ is the space of all sequences of scalars $x=(a_1,
a_2, \ldots)$ such that $\sum_{n=1}^{\infty}M(|a_n|/\rho)<\infty$
for some $\rho>0$, equipped with the norm
$$\|x\|=\inf \big\{\rho>0 :  \sum_{n=1}^{\infty}M(|a_n|/\rho)\le
1\big\}.$$ We will always assume that $M$ satisfies the
\textit{$\Delta_2$-condition} at zero (i.e., that there exists $C>0$
such that $M(2t) \le CM(t)$ for all $0 \le t \le 1/2$). Then the
unit vectors form a normalized symmetric basis for $\ell_M$. If $N$
also satisfies the $\Delta_2$-condition at zero then $M$ and $N$ are
\textit{equivalent} if there exists a constant $C>0$ such that
$(1/C) N(t) \le M(t) \le CN(t)$ for all $0\le t \le 1$.

If $C\ge1$ and $M$ and $N$ are two Orlicz functions such that $N(t)
\le CM(t)$ for all $0<t\le1$, then the unit vector basis of $\ell_M$
$C$-dominates that of $\ell_N$. Conversely, if $M$ and $N$ satisfy
the $\Delta_2$-condition at zero and the unit vector basis of
$\ell_M$ dominates that of $\ell_N$ then there exists $C\ge1$ such
that $N(t) \le CM(t)$ for all $0<t\le1$.

If $M$ satisfies the $\Delta_2$-condition at zero  then an Orlicz
sequence space $\ell_N$ is isomorphic to a subspace of $\ell_M$ if
and only if $N$ is equivalent to some function  in $C_{M, 1}$, where
$C_{M,1}$ is the norm-closed convex hull in $C[0,1]$) of the set
\begin{equation} \label{eq: EM1}
E_{M,1}=\overline {\Big\{\frac{M(\lambda t)}{M(\lambda)};\ 0<\lambda
<1\Big\}}. \end{equation} See \cite[Lemma 4.a.6 and remark (p.
141)]{LT} for this result.

As noted in \cite{S1}, this is easily generalized to the spreading
models of $\ell_M$:  $(\tilde {x_i})$ is a spreading model generated
by a normalized block sequence $(x_i)$ in $\ell_M$ if and only if
$(\tilde {x_i})$ is \textit{isometrically} equivalent to the unit
vector basis of $\ell_N$ for some $N\in C_{M,1}$.

We will use the following method of representing Orlicz functions by
sequences of zeros and ones, introduced by Lindenstrauss and
Tzafriri \cite[p. 161]{LT}.

Fix $0<\tau<1$ and $1<r<p<\infty$. For every sequence of zeros and
ones, $\eta=(\eta(n))_{n=1}^{\infty}$ (i.e. $\eta(n)\in\{0,1\}$ for
all $n$), let $M_{\eta}$ be the  piecewise linear function defined
on $[0,1]$ satisfying $M_{\eta}(0)=0$, $M_{\eta}(1)=1$, and
$$M_{\eta}(\tau^k)=\tau^{r
  k+(p-r)\sum_{n=1}^k\eta(n)}, \qquad k=1,2,\ldots.$$
\begin{lem} \label{lem: prlemma}
Fix $0<\tau<1$ and $1<r<p<\infty$. If $p-r$ is sufficiently small,
then, for all $\eta$, $M_\eta$ is an Orlicz function satisfying the
$\Delta_2$-condition at zero. \end{lem}
\begin{proof} To show that $M_\eta$ is convex it suffices to check
that the slope of the chord joining
$(\tau^{n+1},M_\eta(\tau^{n+1}))$ to $(\tau^n,M_\eta(\tau^n))$ is a
decreasing function of $n$, i.e.
\begin{equation} \label{eq: chordslope}
\frac{M_\eta(\tau^{n})-M_\eta(\tau^{n+1})}{\tau^n-\tau^{n+1}} \le
\frac{M_\eta(\tau^{n-1})-M_\eta(\tau^{n})}{\tau^{n-1}-\tau^{n}}
\end{equation}
Using the fact that $M_\eta(\tau^{k+1})= \tau^r M_\eta(\tau^{k})$ if
$\eta(k+1)=0$ and $M_\eta(\tau^{k+1})= \tau^p M_\eta(\tau^{k})$ if
$\eta(k+1)=1$, \eqref{eq: chordslope} simplifies to the following
pair of conditions: \begin{equation} \label{eq: chordslope2}
\tau^{r-1}(1-\tau^p) \le 1-\tau^r \quad\text{and}\quad
\tau^{p-1}(1-\tau^r)\le 1-\tau^p. \end{equation} Both conditions are
clearly satisfied if $p-r$ is sufficiently small. (Note that the
first condition is \textit{not} satisfied, however, if $r$ is very
close to  $1$.) The $\Delta_2$-condition is easily checked.
\end{proof}
Henceforth, we shall always assume that $\tau$, $p$, and $r$
satisfy the conclusion of
Lemma~\ref{lem: prlemma}.

\begin{prop} \label{prop: L1embedding} Suppose that $1< p <2$.
Then there exists $C<\infty$ such that for each sequence $\eta$ of
zeros and ones, $\ell_{M_\eta}$ $C$-embeds into $L_1[0,1]$.
\end{prop}
\begin{proof} Observe that the inequalities in \eqref{eq:
chordslope2} are reversed if $0<r<p<1$. This implies that
$M_\eta(\sqrt t)$ is equivalent to a {\it concave} function if
$p<2$. By  a result of Bretagnolle and Dacuhna-Castelle \cite{BD}
$\ell_{M_\eta}$ embeds isomorphically into $L_1[0,1]$. To see that
there is a uniform embedding constant, observe that there exists a
`universal' sequence $\rho = (\rho(n))_{n=1}^\infty$ such that every
sequence of zeros and ones, $\eta$, is a pointwise limit of the
collection $\{(\rho(n+k)_{n=1}^\infty \colon k \in \mathbb{N}\}$ of
left shifts of $\rho$: indeed, let $\rho$  be the concatenation of
all possible finite sequences of zeros and ones. It follows that
$E_{M_\rho,1}$ (see \eqref{eq: EM1}) contains $M_\eta$ for every
$\eta$, and hence that $\ell_{M_\eta}$ is isometric to a spreading
model of $\ell_{M_\rho}$  for every $\eta$.  Finally, since
$\ell_{M_\rho}$ $C$-embeds into $L_1[0,1]$ for some $C<\infty$, it
follows that $\ell_{M_\eta}$ $C$-embeds into $L_1[0,1]$ for every
$\eta$.
\end{proof}

We will be interested in only a simple class
of such spaces as described in the following.

\begin{lem}\label{3.2}
Let $1<r<p<\infty$, and $0<\tau<1$. Let $(n_k)\subset\nat$ satisfy
$n_1=1$ and $n_{k+1}-n_k\uparrow\infty$, and put
\begin{equation*}\rho(i)  = \begin{cases} 0
&\text{if $i= n_k$}\\
1 &\text{otherwise.} \end{cases}
\end{equation*}
Let $M:=M_{\rho}$ be the corresponding Orlicz function. Then
$\ell_M$ satisfies the following: \begin{itemize}
\item[(a)] Every spreading model is $\tau^{-2p}$-dominated by the
unit vector basis of $\ell_M$.
\item[(b)] Every spreading model is equivalent either to the unit
vector basis of $\ell_M$ or to the unit vector basis of $\ell_p$.
\item[(c)] Every spreading model that is equivalent to the unit
vector basis of $\ell_p$ is actually $\tau^{-5p}$-equivalent to
the unit vector basis of $\ell_p$.
\item[(d)] $\ell_M$ is reflexive.
\end{itemize}
\end{lem}

\begin{proof}
Observe that
\begin{equation*}\sum_{i=1}^{n} \rho(i)\le
\sum_{i=k+1}^{k+n}\rho(i)\ \ {\rm for\ all}\ k,
n\in\nat.\end{equation*} Therefore, for all $\lambda=\tau^k$ and
$t=\tau^n$, we have
\begin{eqnarray*}
\frac{M(\lambda t)}{M(\lambda)}&=&
\frac{\tau^{r(k+n)+(p-r)\sum_{i=1}^{k+n}\rho(i)}}
{\tau^{rk+(p-r)\sum_{i=1}^{k}\rho(i)}}\nonumber\\
&=&\tau^{rn+(p-r)\sum_{i=k+1}^{k+n}\rho(i)}
\le\tau^{rn+(p-r)\sum_{i=1}^{n}\rho(i)}=M(t).
\end{eqnarray*}
A simple calculation now yields, for all
$0<\lambda, t<1$, that
\begin{equation}\label{ineq1}
\frac{M(\lambda t)}{M(\lambda)}\le\tau^{-2p}M(t).
\end{equation}

Now let $N\in C_{M,1}$. Then, by the definition of $C_{M,1}$, $N$ is
the limit in the uniform norm of a sequence of convex combinations
$(F_n)$ of the form
$$F_n=\sum_{i\in A_n}a_i\frac{M(\lambda_i t)}{M(\lambda_i)},$$
for some finite $A_n\subset\nat$, $0<\lambda_i\le 1$,
and positive $(a_i)$ with $\sum_{i\in A_n}a_i=1$.

Thus (\ref{ineq1})  implies that for all $N\in
C_{M,1}$,
\begin{equation}\label{uniform-ineq}
N(t)\le \tau^{-2p}M(t)\ \ {\rm for\ all}\ 0<t<1,
\end{equation}
which proves (a).
To see that $N$ is equivalent  either to $M$ or to $t^p$, we
distinguish two cases corresponding to
the manner in which the  sequence $(F_n)$  converges to $N$.

For the first case, suppose that
there exists $n_0\in\nat$, $\bar\lambda>0$  and $\delta>0$ such
that for all $n\ge n_0$,
$$\sum_{\lambda_i\ge \bar\lambda, i\in A_n}a_i\ge\delta.$$
It follows that
$$N(t)\ge \delta M(\bar\lambda
  t),\ {\rm for\ all}\ t>0,$$
which, along with (\ref{uniform-ineq}), implies that $N$ is
equivalent to $M$.

For the second case, we
suppose  that for all $n_0\in\nat$, $\bar\lambda>0$ and $\delta>0$,
there exists $n \ge n_0$ such that
$$\sum_{\lambda_i\ge\bar\lambda, i\in A_n}a_i<\delta.$$
Fix $t = \tau^m$. Since $n_{k+1}-n_k \uparrow \infty$, it follows
that every $m$ consecutive terms of $\rho$ which begin sufficiently
far along the sequence  can contain at most one zero term. This
implies that if $\lambda=\tau^k$ is sufficiently small then
\begin{equation} \label{eq: lpcase} t^p \le  \frac{M(\lambda
t)}{M(\lambda)} \le \tau^{-p} t^p. \end{equation} Now fix $\delta>0$
and $0<t<1$. It follows easily from \eqref{eq: lpcase} that there
exists $\bar\lambda>0$ such that for all $\lambda < \bar\lambda$
\begin{equation} \label{eq: lpcase2}
\tau^{2p}t^p \le  \frac{M(\lambda t)}{M(\lambda)} \le \tau^{-3p} t^p.
\end{equation} By assumption there exists $n \in \mathbb{N}$ such that
\begin{equation} \label{eq: FandN}
|F_n(t)-N(t)| < \delta \quad\text{and}\quad
\sum_{\lambda_i\ge\bar\lambda, i\in A_n}a_i<\delta. \end{equation}
Now \eqref{eq: lpcase2} gives \begin{equation*} F_n(t)
=\sum_{\lambda_i<\bar\lambda, i\in A_n}a_i \frac{M(\lambda_i
t)}{M(\lambda_i)} +\sum_{\lambda_i\ge\bar\lambda, i\in A_n}a_i
\frac{M(\lambda_i t)}{M(\lambda_i)} \le \tau^{-3p} t^p + \delta.
\end{equation*} A similar calculation yields
$$F_n(t) \ge (1-\delta) \tau^{2p} t^p - \delta.$$
Since  $|N(t)-F_n(t)| < \delta$ and $\delta>0$ is arbitrary, we get
$$ \tau^{2p} t^p \le N(t) \le \tau^{-3p}t^p.$$
This proves (c) and also completes the proof of (b). Finally, (b)
and \cite[Lemma 4.a.6]{LT} imply that $\ell_1$ is not isomorphic to
a subspace of $\ell_M$, which in turn implies by \cite[Proposition
4.a.4]{LT} that $\ell_M$ is reflexive, which proves (d).
\end{proof}

We will also make use of the following general fact.

\begin{lem}\label{3.3}
Let $X=\big(\sum_{j=1}^{\infty}\oplus X_j\big)_p$, where $1\le p
<\infty$ and each $X_j$ is an infinite-dimensional Banach space, and
let $(\tilde x_i)\in SP_w(X)$ be a spreading model generated by a
normalized weakly null sequence in $X$.  Then there exist
non-negative $(c_j)_{j=0}^{\infty}$ with
$\sum_{j=0}^{\infty}c_j^p=1$ and normalized spreading models
$(\tilde x_i^j)_i \in SP_w(X_j)$ such that for all scalars $(a_i)$
$$\big\|\sum_i a_i \tilde
x_i\big\|=\Big[\sum_{j=1}^{\infty}c_j^p\big\|\sum_i a_i\tilde
x_i^j\big\|^p + c_0^p\sum_i |a_i|^p\Big]^{1/p}.$$
\end{lem}

\begin{proof}
Suppose that the normalized weakly null sequence $(y_i)$ generates
the spreading model $(\tilde x_i)$. Write
$y_i=(y_i^j)_{j=1}^{\infty}$, where $y^j_i \in X_j$ for each $j$. By
a diagonalization argument, $(y_i)$ has a subsequence $(x_i)$ such
that
\begin{equation}\label{cantor1}
\lim_{i\to\infty}\|x_i^j\|=c_j\ \ {\rm and}\
\sup_{i\ge
j}\big|\|x_i^j\|-c_j\big|\le \frac{1}{2^j},
\end{equation}
and $(c_j^{-1}x_i^j)_{i=1}^{\infty}$ generates a normalized
spreading model $(\tilde x_i^j)_{i} \in SP_w(X_j)$. Note that
(\ref{cantor1})  implies that $c_0:=\lim_{i\to
\infty}\|x_i-P_i(x_i)\|$ exists, where $P_i(x_i)=(x_i^1, x_i^2,
\ldots, x_i^i, 0, 0, 0, \ldots)$. One now checks that the spreading
model $(\tilde x_i)$ generated by $(x_i)$ is given by the stated
formula. (Note also that $\sum_{j=0}^{\infty}c_j^p=1$ since $(x_i)$
is normalized.)
\end{proof}

Before proving Theorem~\ref{thm: lattice} we give an application of
Lemma~\ref{3.3}.

\begin{thm} \label{prop: ordinal}
Let $X$ be an infinite-dimensional Banach space such that
$SP_w(X)$ is a countable chain. Then there exists a countable
ordinal $\alpha$ such that $SP_w(X)$ is order-isomorphic to
$\alpha$ with the \emph{reverse order}. Conversely, if
$\alpha\ge1$ is a countable ordinal then there exists a reflexive
Banach space $X$ such that $SP_w(X)$ is order-isomorphic to
$\alpha$ with the reverse order. \end{thm}
\begin{proof}
For the first part,  by Fact~\ref{fact1.3} $SP_w(X)$ does not admit
an infinite strictly increasing sequence. Thus the reverse order on
$SP_w(X)$ is a well-ordering and hence is order-isomorphic to a
countable ordinal.  For the converse, let $\beta \mapsto p_\beta$
($\beta < \alpha$) be an increasing order-isomorphism from $\alpha$
onto a subset of $[2,3]$ such that $p_0=2$, and set $X :=
(\sum_{\beta<\alpha} \oplus \ell_{p_\beta})_2$. Using the
well-foundedness of $\alpha$ and the monotonicity property of  the
$\ell_p$ norms (i.e.,  that $\|\cdot\|_q \le \|\cdot\|_p$ if $p \le
q$), it follows easily from   Lemma~\ref{3.3} that every normalized
spreading model of $X$ is equivalent to the unit vector basis of
$\ell_{p_\beta}$ for some $\beta<\alpha$; so $SP_w(X)$ is
order-isomorphic to $\alpha$ with the reverse order.
\end{proof}

We now proceed to the

\begin{proof}[Proof of Theorem 3.1]
For  convenience we shall assume that $L$ is countably infinite.
(When $L$ is finite only minor notational changes are needed.)  The
space $X_L$ will be of the form
$X_L=\left(\sum_{j=0}^{\infty}\oplus\ell_{M_j}\right)_p$ for
suitably constructed Orlicz sequence space $\ell_{M_j}$'s, with $M_j
:= M_{\rho_j}$ for certain sequences $\rho_j$ of zeros and ones (for
the same $\tau, r$, and $p$). The `patterns' of the $\rho_j$'s will
be of the form
\begin{equation*}\rho_j(i)  = \begin{cases} 0
&\text{if $i \in \sigma(j)$}\\
1 &\text{otherwise.} \end{cases}
\end{equation*}
for some fast increasing sequence $\sigma(j)\subset\nat$, with $1
\in \sigma(j)$. For simplicity, for every $j$ we will take
$\sigma(j)$ to be a subset of $\mathcal{M}=\{1,2, 2^2, 2^3,\ldots
\}$ which will ensure that the hypothesis of Lemma~\ref{3.2} is
satisfied.

The  patterns of the $\rho_j$'s (equivalently, the $\sigma_j$'s)
will be developed inductively on finite intervals of $\mathbb{N}$
according to a two-step  procedure which we call
\textit{$(\varepsilon,A)$-domination}.

Let $A\subset\nat$ and $\ep>0$. Suppose that for some
$N\in\mathbb{N}$,  the $\rho_j$'s have already been defined on
the
initial segment $[1, N]$ so that
\begin{equation}\label{balanced}
\sum_{i=1}^N \rho_j(i)=\sum_{i=1}^N \rho_k(i),\ {\rm for\ all}\ j,
k\in\nat.
\end{equation}

The $(\varepsilon,A)$-domination procedure extends the definition of
the $\rho_j$'s to an intial segment $[1,N_1]$ for some $N_1>N$. Let
us first dispose of some trivial cases. If $A = \emptyset$ or if $A
= \mathbb{N}$ then set $N_1=N+1$ and $\rho_j(N_1)=1$ for all $j$.

Now suppose that both $A$ and $\mathbb{N} \setminus A$ are
non-empty. The first step of the procedure is carried out  as
follows. Choose a sufficiently large (just how large is specified
below) integer $m>N$. For all $k\in \nat\setminus A$ place $0$'s
on the coordinates from $[N+1, m]\cap\mathcal{M}$  of the
$\rho_k$'s (while  the rest of the coordinates of the interval are
filled with $1$'s), and for all $j\in A$ place $1$'s on \emph{all}
the coordinates from $[N,m]$ of  the $\rho_j$'s , where $m$ is
chosen so that
\begin{equation*}
\sum_{i=1}^m \rho_j(i)-\sum_{i=1}^m \rho_k(i)
\end{equation*}
is sufficiently large to ensure that
\begin{equation*}
\frac{M_{\rho_j}(\tau^m)}{M_{\rho_k}(\tau^m)}<\ep,\ {\rm for\ all}\
j\in A, \ k\in\nat\setminus A.
\end{equation*}
For the second step we choose a sufficiently large integer $N_1>m$
(just how large is specified below), with $N_1 \in \mathcal{M}$, and
`rebalance' all of the $\rho_j$'s on the interval $[m+1,N_1]$. This
is achieved by placing $0$'s on the coordinates from $[m+1,N_1] \cap
\mathcal{M}$ for all the $\rho_j$'s ($j \in A$) and by placing $1$'s
on the coordinates from $[m+1,N_1]$ for all the $\rho_k$'s ($k \in
\nat\setminus A$), where $N_1$ is chosen so that \eqref{balanced} is
satisfied with $N$ replaced by $N_1$. At the end of this second step
the $M_j$'s are equal again, i.e.
$$M_j(\tau^{N_1}) = M_k(\tau^{N_1}) \quad \text{for all $j.k \in
\mathbb{N}$}$$

We now pass to the main construction.
Let $L=\{e_0, e_1, e_2, \ldots,\}$ be the given countable lattice,
where $e_0$ is the minimum element. Consider $\bar L=\{\bar e_1, \bar
e_2, \ldots\}$, where $\bar e_j=\{i\in\nat : e_i\le e_j\}$ for all
$j\in\nat$. Put $\rho_0=(1, 1, 1, \dots)$.

We begin by setting $\rho_j(1)=0$ for all $j\in \mathbb{N}$, which
ensures that the $\rho_j's$ satisfy Lemma~\ref{3.2}. Now, for every
$j\in\nat$ and every $\ep=2^{-k}$, $k=1, 2, \dots$, we carry out an
$(\ep, A)$-domination procedure for $A=\bar e_j$. Since there are
countably many choices we can enumerate some order in which to carry
out all $(\ep, A)$-dominations.

The resulting sequences $\rho_0, \rho_1, \rho_2, \ldots$ have the
following properties.

(i) $M_{\rho_0}$ is equivalent to the function $t^p$.

(ii) For all $i,j \in \mathbb{N}\cup\{0\}$,
there exists a constant $C<\infty$ such that
$$M_{\rho_i}(t) \le CM_{\rho_j}(t) \quad\text{for all $0<t<1$}$$
if and only if $e_i \le e_j$ in $(L,\le)$. Moreover, if there exists
such a $C$ then $C=1$ works.

(iii) For every non-empty finite set $F\subset \nat \cup \{0\}$
\begin{equation*}
\max_{j\in F} M_{\rho_j}= M_{\rho_{j_0}}, \ {\rm where}\
e_{j_0}=\bigvee_{j\in F}e_j.
\end{equation*}
\noindent{\emph{Proof of (iii).}} To derive a contradiction, assume
that there exists $t=\tau^m$ such that $\max_{j\in F} M_{\rho_j}(t)<
M_{\rho_{j_0}}(t)$. Because of the `rebalancing' step in the
domination procedure, it follows that  $m$ belongs to an interval of
$\mathbb{N}$ where an $(\ep, A)$-domination takes place for some $A$
such that $F\subseteq  A$ and $j_0\in\nat\setminus A$. There exists
$k \in \mathbb{N}$ such that $A= \bar e_k$. Then $e_j \le e_k$ for
all  $j \in F$. Since $L$ is a lattice it follows that $e_{j_0} \le
e_k$, and hence $j_0 \in A$, which is the desired contradiction.

(iv) Let $B$ be a non-empty subset of  $\nat\cup\{0\}$  Then there
exists a finite subset $F$ of $B$ such that
$$\max_{j\in B} M_{\rho_j}=\max_{j\in F} M_{\rho_j}.$$

\noindent{\emph{Proof of  (iv).}} Suppose not. Then there exists
  $(j_k)_{k=1}^{\infty}\subset B$ such that for all $n \in \mathbb{N}$
$$\max_{1\le k\le n}M_{\rho_{j_k}} < \max_{1\le k\le n+1} M_{\rho_{j_k}}.$$
This, however, implies by (iii) that
$$\bigvee_{1\le k\le n} e_{j_k}<\bigvee_{1\le j\le n+1}e_{j_k},\ {\rm for\
each}\ n.$$ But this contradicts our assumption that there are no
increasing infinite sequences in $L$. \qed

Now consider
$$X_L=\Big(\sum_{j=0}^{\infty}\oplus\ell_{M_j}\Big)_p,$$
where $M_j=M_{\rho_j}$, $j\in\nat\cup\{0\}$. By (d) of
Lemma~\ref{3.2} each $\ell_{M_j}$ is reflexive and hence $X_L$ is
also reflexive. It follows from property (ii) that the
\emph{collection} of spreading models generated by the unit vector
bases  of each $\ell_{M_j}$ is order-isomorphic to $L$. Therefore it
remains to show that \textit{every} spreading model of $X_L$ is
equivalent to the unit vector basis $(b_i^j)_i$ of $\ell_{M_j}$ for
some $j\in\nat\cup\{0\}$.

Let $(\tilde x^i)$ be a normalized spreading model of $X_L$. Then,
for all $(a_i)\in c_{00}$, we have by Lemma~\ref{3.3}
\begin{equation} \label{eq: spmodelformula}
\big\|\sum_i a_i \tilde x_i\big\|= \Big[\sum_j c_j^p \big\|\sum_i
a_i \tilde x_i^j\big\|^p + c_0^p \sum_i |a_i|^p\Big]^{1/p},
\end{equation}
where $(\tilde x_i^j)_i$ is a normalized spreading model of
$\ell_{M_j}$ and $(c_j)_{j=0}^\infty$ belongs to the non-negative
unit sphere of $\ell_p$.

Let $B$ be the collection of all $j \in \mathbb{N}$ such that $c_j
\ne 0$ \textit{and} such that $(\tilde x_i^j)_i$ is equivalent to
$(b_i^j)_i$. If $j \notin B$ then either $c_j=0$ or, by
Lemma~\ref{3.2}, $(\tilde x_i^j)_i$ is $\tau^{-5p}$-equivalent to
the unit vector basis of $\ell_p$. Thus, if $B = \emptyset$, then
\eqref{eq: spmodelformula} implies that $(\tilde x^i)$ is equivalent
to the unit vector basis of $\ell_p$ and hence  to $(b_i^0)_i$. So
suppose that $B \ne \emptyset$. Then, by Lemma~\ref{3.2}, each
$(\tilde x_i^j)_i$ ($j \in B$) is $\tau^{-2p}$-dominated by
$(b_i^j)_i$. By properties  (iii) and (iv) above there exist a
finite set $F \subset B$ and $j_0 \in \mathbb{N}$ such that
$$\max_{j\in B} M_{\rho_j}=\max_{j\in F} M_{\rho_j}=M_{\rho_{j_0}}.$$
Hence there exists $0 \le K <\infty$ such that \begin{align*}
\big\|\sum_i a_i \tilde x_i\big\| &\le \Big[\tau^{-2p}\sum_{j\in B}
c_j^p \big\|\sum_i
  a_i b_{i}^j\big\|^p + K \sum_i |a_i|^p\Big]^{1/p}\\
&\le \Big[(\tau^{-2p}\sum_{j\in B} c_j^p) \big\|\sum_i
  a_i b_{i}^{j_0}\big\|^p + K \sum_i |a_i|^p\Big]^{1/p},
\end{align*}
which implies that $(\tilde x_i)_i \le (b_i^{j_0})$. On the other
hand, there exists $c>0$ such that
\begin{align*}
\big\|\sum_i a_i \tilde x_i\big\| &\ge \Big[\sum_{j\in F} c_j^p
\big\|\sum_i a_i \tilde x_i^j\big\|^p\Big]^{1/p}\\
& \ge c \max_{j \in F} \big\|\sum_i a_i b_i^j\big\|\\
\intertext{(since  $(\tilde x_i^j)_i$ is equivalent to $(b^j_i)_i$
for each $j \in F$)} &\ge \frac{c}{\operatorname{card} F}
\big\|\sum_i a_i b_i^{j_0}\big\|. \end{align*} Thus,  $(\tilde
x_i)_i$ is equivalent to $(b_i^{j_0})_i$ .
\end{proof}
\begin{remark} For each $1<p<\infty$, the above construction
yields the unit vector basis of  $\ell_p$ as the minimum element of
$SP_w(X_L)$. If we allow $X_L$ to be nonreflexive we can obtain
$c_0$ as the minimum element. However, this requires a rather
different construction. Using  results of Casazza and Lin \cite{CL},
it is possible to construct a $c_0$-sum of duals of certain Lorentz
sequence spaces for which $c_0$ is the minimum element of
$SP_w(X_L)$. We omit the details of this result. \end{remark}
\begin{remark} Let $\rho$ be the universal sequence used in the proof of
Proposition~\ref{prop: L1embedding}. It follows  from the proof of
Theorem~\ref{thm: lattice} that $SP_w(\ell_{M_\rho})$ contains a
\textit{subset} that is order-isomorphic to any given countable
poset $P$. (Note also that there is a universal countable poset.)
By Proposition~\ref{prop: L1embedding},  $\ell_{M_\rho}$ is
isomorphic to a reflexive subspace of $L_1[0,1]$ when $p<2$.
 \end{remark}
\begin{cor} For every finite lattice $L$ there exists a reflexive
space $X_L$ such that $SP_w(X_L)$ is order-isomorphic to $L$.
\end{cor}

\begin{cor} Let $L$ be a finite lattice (resp.\ countable lattice with
 a minimum element and without any infinite increasing sequence).
There exists a reflexive (resp.\ non-reflexive) subspace $Y_L$ of $L_1[0,1]$
such that $SP_w(Y_L)$ is order-isomorphic to $L$.
\end{cor}
\begin{proof}  Using the notation of Theorem~\ref{thm: lattice}, let
$$Y_L=\Big(\sum_{j=0}^{\infty}\oplus\ell_{M_j}\Big)_1.$$
By Proposition~\ref{prop: L1embedding}, if $p<2$  then for some
$C<\infty$ each $\ell_{M_j}$ $C$-embeds into $L_1[0,1]$, and hence
$Y_L$ is isomorphic to a subspace of $L_1[0,1]$. Moreover, $Y_L$
is reflexive if and only if $L$ is finite. The proof of
Theorem~\ref{thm: lattice} shows that $SP_w(Y_L)$ is
order-isomorphic to $L$. (Note that if $L$ is infinite then $Y_L$
also has an  $\ell_1$ spreading model that  is \textit{not}
generated by a weakly null sequence.)
\end{proof}


\end{document}